\documentclass[10pt,leqno]{amsart}
\usepackage{indentfirst, amssymb,amsmath,amsthm}
\newtheorem{theo}{Theorem}[section]
\newtheorem{lem}{Lemma}[section]
\newtheorem{ques}{Question}[section]
\newtheorem{cor}{Corollary}[section]

\newtheorem{exm}{Example}[section]
\newtheorem{defi}{Definition}[section]
\newtheorem{rem}{Remark}[section]
\newcommand{\ol}{\overline}
\newcommand{\be}{\begin{equation}}
\newcommand{\ee}{\end{equation}}
\newcommand{\beas}{\begin{eqnarray*}}
\newcommand{\eeas}{\end{eqnarray*}}
\newcommand{\bea}{\begin{eqnarray}}
\newcommand{\eea}{\end{eqnarray}}

\numberwithin{equation}{section}
\begin{document}
\title[ Uniqueness of meromorphic functions concerning derivatives]{Uniqueness of Meromorphic functions concerning k-th derivatives and difference operators}
\date{}
\author[G. Haldar ]{ Goutam Haldar}
\date{}
\maketitle
\let\thefootnote\relax
\footnotetext{Department of Mathematics, Malda College, Rabindra Avenue, Malda, West Bengal 732101, India}
\footnotetext{\textbf{E-mail:} goutamiit1986@gmail.com}
\setcounter{footnote}{0}


\begin{abstract}
	In this paper, we continue to study the sharing value problems for higher order derivatives of meromorphic functions with its linear difference and $q$-difference operators. Some of our results generalize and improve the results of Meng--Liu (J. Appl. Math. and Informatics, 37(2019), 133--148) to a large extent. 
	
	\vskip 0.4 true cm

	\noindent
	
	\noindent
	
\end{abstract}
\section{Introduction}
Let $f$ and $g$ be two non-constant meromorphic functions defined in the open complex plane $\mathbb{C}$. If for some $a\in\mathbb{C}\cup\{\infty\}$, $f-a$ and $g-a$ have the same set of zeros with the same multiplicities, we say that $f$ and $g$ share the value $a$ CM (counting multiplicities), and if we do not consider the multiplicities then $f$ and $g$ are said to share the value $a$ IM (ignoring multiplicities). We assume that the readers are familiar with the standard notations symbols such as $T(r,f)$, $N(r,a;f)$ ($\overline N(r,a;f)$) of Nevanlinna’s value distribution theory (see \cite{Hayman & 1964}).\vspace{1mm}
\par In 2001, Lahiri (\cite{Lahiri & Complex Var & 2001}, \cite{Lahiri & Nagaya & 2001}) introduced the definition ofweighted sharing, which plays a key role in uniqueness theory as far as relaxation of sharing is concerned. In the following we explain the notion.
\begin{defi}\cite{Lahiri & Complex Var & 2001}
	Let $k$ be a non-negative integer or infinity. For $a\in \mathbb{C}\cup\{\infty\}$ we denote by $E_{k}(a,f)$ the set of all $a$-points of $f$, where an $a$ point of multiplicity $m$ is counted $m$ times if $m\leq k$ and $k+1$ times if $m>k.$ If $E_{k}(a,f)=E_{k}(a,g),$ we say that $f$, $g$ share the value $a$ with weight $k$. 
\end{defi}
\par We write $f$, $g$ share $(a,k)$ to mean that $f,$ $g$ share the value $a$ with weight $k.$ Clearly if $f,$ $g$ share $(a,k)$ then $f,$ $g$ share $(a,p)$ for any integer $p$, $0\leq p<k.$ Also we note that $f,$ $g$ share a value $a$ IM or CM if and only if $f,$ $g$ share $(a,0)$ or $(a,\infty)$ respectively.\vspace{1mm}
\begin{defi}\cite{Lahiri & 2001}
	For $a\in \mathbb{C}\cup \{\infty\},$ we denote by $N(r,a;f\mid=1)$ the counting function of simple $a$-points of $f.$ For a positive integer $m,$ we denote by $N(r,a;f\mid\leq m)$ $(N(r,a;f\mid\geq m))$ the counting function of those $a$-point of $f$ whose multiplicities are not greater (less) than $m$, where each $a$-point is counted according to its multiplicity.
\end{defi}
$\overline N(r,a;f\mid\leq m)$ $(\overline N(r,a;f\mid\geq m))$ are defined similarly except that in counting the $a$-points of $f$ we ignore the multiplicity. Also $N(r,a;f\mid< m)$, $N(r,a;f\mid> m),$ $\overline N(r,a;f\mid< m)$ and $\overline N(r,a;f\mid>m)$ are defined similarly.
\begin{defi}\cite{Lahiri & Complex Var & 2001}
	We denote by $N_2(r,a;f)$ the sum $\overline N(r,a;f)+\overline N(r,a;f\mid\geq 2)$.
\end{defi}
\begin{defi}\cite{Lahiri & Complex Var & 2001}
	Let $f$ and $g$ share a value $a$ IM. We denote by $\overline N_{*}(r,a;f,g)$ the counting function of those $a$-points of $f$ whose multiplicities differ from the multiplicities of the corresponding $a$-points of $g$.\end{defi}
\par Let $c$ be a nonzero complex constant, and let $f(z)$ be a meromorphic function. The shift operator is denoted by $f(z+c)$. Also, we use the notations $\Delta_cf$ and $\Delta_c^kf$ to denote the
difference and kth-order difference operators of $f$, which are respectively defined as \begin{eqnarray*} \Delta_cf=f(z+c)-f(z),\;\; \Delta_c^kf(z)=\Delta_c(\Delta_c^{k-1}f(z)),\;\; k\in \mathbb{N},\;k\geq2.\end{eqnarray*} We note that $\Delta_cf$ and $\Delta_c^kf$ are nothing but linear combination of different shift operators. So for generalization of those operators, it is reasonable to introduce the linear difference operators $L(z,f)$ as follows:
\begin{eqnarray} \label{e1.1} L(z,f)=\sum_{j=0}^{p}a_jf(z+c_j),\end{eqnarray} where $p\in \mathbb{N}\cup\{0\}$ and $a_j$ and $c_j$'s are complex constants with at-least one $a_j$'s are non-zero.\vspace{1mm}
\par For a non-zero complex constant $q$ and a meromorphic function $f$, the $q$-shift and $q$-difference operators are defined, respectively by $f(qz)$ and $\Delta_qf=f(qz)-f(z)$. Here also we generalize these operators as follows:
\begin{eqnarray}\label{e1.2}  L_q(z,f)=\sum_{j=0}^{r}b_jf(q_jz+d_j),\end{eqnarray} where $r$ is a non-negative integer, and $q_j$, $ b_j$, $d_j$'s are complex constants with at-least one of $b_j$  is non-zero.\vspace{1mm}
\par It was Rubel--Yang \cite{Rubel & Yang & 1977} who first initiated the problem of uniqueness of meromorphic functions sharing two values, and obtained the following result.
\begin{theo}\label{t1.5}\cite{Rubel & Yang & 1977}
	Let $f$ be a non-constant entire function. If $f$ shares two distinct finite values CM with $f^{\prime}$, then $f\equiv f^{\prime}$.	
\end{theo}
\par Mues--Steinmetz \cite{Mues & Steinmetz & 1979} improved the above result by relaxing the nature of sharing two values from CM to IM.  After that Mues--Steinmetz \cite{Mues & Steinmetz & 1983}, and Gundersen \cite{Gundersen & JMAA & 1980} improved Theorem A to non-constant meromorphic functions.\vspace{1mm}
\par Recently, the difference analogue of classical Nevanlinna theory for meromorphic functions of finite order was established by Halburd--Korhonen \cite{Halburd & Korhonen & Ann Acad & 2006, Halburd & Korhonen & JMAA & 2006}, Chiang--Feng \cite{Chiang & Feng & 2008}, independently, and developed
by Halburd--Korhonen--Tohge \cite{ Halburd & Korhonen & Tohge & 2014} for hyper order strictly less than 1. After that, there has been an increasing interest in studying the uniqueness problems of meromorphic functions related to their shift or difference operators (see \cite{Chen & Sci. Chaina & 2014, Halburd & Korhonen & 2017, Laine & Berlin & 1993, Liu & Yang 2013, Qi & Yang & 2020, Banerjee & Bhattacharya & 2019, Bhusnurmath & Kabbur & 2013, Chen & Chen & 2012, Li & Chen & 2014, Zhang & JMAA & 2010}).\vspace{1mm}
\par As we know that the time-delay differential equation $f(x)=f(x-k)$, $k>0$ plays an important roll in real analysis, and it has been rigorously studied. For complex variable counterpart, Liu-Dong \cite{Liu & Dong & 2014} studied the complex differential-difference equation $f(z)=f(z+c)$, where $c$ is a non-zero constant.\vspace{1mm}
\par In 2018, Qi \textit{et al.} \cite{Qi & Li & Yang & 2018} looked at this complex differential-difference equation from a different perspective. In fact, they considered the value sharing problem related to $f^{\prime}(z)$ and $f(z+c)$, where $c$ is a complex number, and obtained the following result.
\begin{theo}\label{t1.6}\cite{Qi & Li & Yang & 2018}
	Let $f$ be a non-constant meromorphic function of finite order, $n\geq9$ be an integer. If $[f^{\prime}(z)]^n$ and $f^n(z + c)$ share $a(\neq0)$ and $\infty$ CM, then $f^{\prime}(z)=tf(z+c)$, for a constant $t$ that satisfies $t^n=1$.
\end{theo}
\par In $2019$, Meng--Liu \cite{Meng & Liu & 2019} reduced the nature of sharing values from CM to finite weight and obtained the following results.
\begin{theo}\label{t1.7}\cite{Meng & Liu & 2019}
	Let $f$ be a non-constant meromorphic function of finite order, $n\geq10$ an integer. If $[f^{\prime}(z)]^n$ and $f^n(z+c)$ share $(1,2)$ and $ (\infty,0)$, then $f^{\prime}(z)=tf(z+c)$, for a constant $t$ that satisfies $t^n =1$.
\end{theo}
\begin{theo}\label{t1.8}\cite{Meng & Liu & 2019}
	Let $f$ be a non-constant meromorphic function of finite order, $n\geq9$ an integer. If $[f^{\prime}(z)]^n$ and $f^n(z+c)$ share $(1,2)$ and $ (\infty,\infty)$, then $f^{\prime}(z)=tf(z+c)$, for a constant $t$ that satisfies $t^n =1$.
\end{theo}
\begin{theo}\label{t1.9}\cite{Meng & Liu & 2019}
	Let $f$ be a non-constant meromorphic function of finite order, $n\geq17$ an integer. If $[f^{\prime}(z)]^n$ and $f^n(z+c)$ share $(1,0)$ and $ (\infty,0)$, then $f^{\prime}(z)=tf(z+c)$, for a constant $t$ that satisfies $t^n =1$.
\end{theo}
For further investigation of Theorems \ref{t1.6}--\ref{t1.9}, we pose the following questions.\vspace{1mm}
\begin{ques}\label{q1}
	Could we determine the relationship between the $k$-th derivative $f^{(k)}(z)$ and the linear difference polynomial $L(z,f)$ as defined in $(\ref{e1.1})$ of a meromorphic (or entire) function $f(z)$ under relax sharing hypothesis?
\end{ques}
\begin{ques}\label{q2}
	Could we further reduce the lower bound of $n$ in Theorems \ref{t1.7}--\ref{t1.9}? 
\end{ques}
In this direction, we prove the following result.
\begin{theo}\label{t1}
	Let $f$ be a non-constant meromorphic function of finite order, $n,\;k$ are positive integers, and $L(z,f)$ are defined in $(\ref{e1.1})$. Suppose  $[f^{(k)}]^n$ and $[L(z,f)]^n$ share $(1,l)$ and $(\infty,m)$, where $0\leq l<\infty$ and $0\leq m\leq\infty$, and one of the following conditions holds:
	\begin{enumerate}
		\item [(i)] $l\geq 2$, $m=0$ and $n\geq 8$;
		\item [(ii)] $l\geq 2$, $m=\infty$ and $n\geq 7$;
		\item [(iii)] $l=1$, $m=0$ and $n\geq 9$;
		\item [(iv)] $l=0$, $m=0$ and $n\geq12$.
	\end{enumerate}
	Then $f^{(k)}(z)=tL(z,f)$, for a non-zero constant $t$ that satisfies $t^n=1$.
\end{theo}
\par We give the following example in the support of Theorem $\ref{t1}$.
\begin{exm}\label{1.1}
	Let $f(z)=e^{z/n}$, where $n$ is a positive integer. Suppose $L(z,f)=f(z+c)+c_0f(z)$, where $c_0$ is a non-zero complex constant such that $c_0\neq1/n$, and $c=n\log((1-c_0n)/n)$.Then one can easily verify that $(f^{\prime})^n$ and $(L(z,f))^n$ satisfy all the conditions of Theorem $\ref{t1}$. Here $f^{\prime}(z)=tL(z,f)$, where $t$ is a constant such that $t^{n}=1$. \end{exm}
\begin{rem}\label{r1}
	Let us suppose that $c_{j}=jc$, $j=0,1,\ldots,p$ and $a_{p}(z)={p\choose 0}$, $a_{p-1}=-{p\choose 1}$, $a_{p-2}={p\choose 2}$. Then from (\ref{e1.1}), it is easily seen that $L(z,f)=\Delta^{p}_{c}f$. Therefore, we obtain the following corollary from Theorem \ref{t1}.
\end{rem} 
\begin{cor}\label{c1}
	Let $f$ be a non-constant meromorphic function of finite order, $n,\;k$ are positive integers, and $L(z,f)$ are defined in $(\ref{e1.1})$. Suppose  $[f^{(k)}]^n$ and $[\Delta^{p}_{c}f]^n$ share $(1,l)$ and $(\infty,m)$, where $0\leq l<\infty$ and $0\leq m\leq\infty$, and one of the following conditions holds:
	\begin{enumerate}
		\item [(i)] $l\geq 2$, $m=0$ and $n\geq 8$;
		\item [(ii)] $l\geq 2$, $m=\infty$ and $n\geq 7$;
		\item [(iii)] $l=1$, $m=0$ and $n\geq 9$;
		\item [(iv)] $l=0$, $m=0$ and $n\geq12$.
	\end{enumerate}
	Then $f^{(k)}(z)=t\Delta^{p}_{c}f$, for a non-zero constant $t$ that satisfies $t^n=1$.
\end{cor}
\par For entire function we prove the following result which is an improvement of Corollary 1.8 of \cite{Meng & Liu & 2019}. 
\begin{theo}\label{t3}
	Let $f$ be a non-constant entire function of finite order, $n,\;k$ are positive integers, and $L(z,f)$ are defined in $(\ref{e1.1})$. Suppose  $[f^{(k)}]^n$ and $[L(z,f)]^n$ share $(1,l)$, and one of the following conditions holds:
	\begin{enumerate}
		\item [(i)] $l\geq 1$ and $n\geq 5$;
		\item [(ii)] $l=0$ and $n\geq 8$;
	\end{enumerate}
	Then $f^{(k)}(z)=tL(z,f)$, for a non-zero constant $t$ that satisfies $t^n=1$.
\end{theo}
\par In the same paper, Meng--Liu \cite{Meng & Liu & 2019} also obtained the following results by replacing $f(z+c)$ with $q$-shift operator $f(qz)$.
\begin{theo}\label{t1.14}\cite{Meng & Liu & 2019}
	Let $f$ be a non-constant meromorphic function of zero order, $n\geq10$ an integer. If $[f^{\prime}(z)]^n$ and $f^n(qz)$ share $(1,2)$ and $ (\infty,0)$, then $f^{\prime}(z)=tf(qz)$, for a constant $t$ that satisfies $t^n =1$.
\end{theo}
\begin{theo}\label{t1.15}\cite{Meng & Liu & 2019}
	Let $f$ be a non-constant meromorphic function of zero order, $n\geq9$ an integer. If $[f^{\prime}(z)]^n$ and $f^n(qz)$ share $(1,2)$ and $ (\infty,\infty)$, then $f^{\prime}(z)=tf(qz)$, for a constant $t$ that satisfies $t^n =1$.
\end{theo}
\begin{theo}\label{t1.16}\cite{Meng & Liu & 2019}
	Let $f$ be a non-constant meromorphic function of zero order, $n\geq17$ an integer. If $[f^{\prime}(z)]^n$ and $f^n(qz)$ share $(1,0)$ and $ (\infty,0)$, then $f^{\prime}(z)=tf(qz)$, for a constant $t$ that satisfies $t^n =1$.
\end{theo}
\par For the generalizations and improvements of Theorems \ref{t1.14}--\ref{t1.16} to a large extent, we obtain the following result.
\begin{theo}\label{t2}
	Let $f$ be a non-constant meromorphic function of zero order, $n,\;k$ are positive integers, and $L_q(z,f)$ are defined in $(\ref{e1.2})$. Suppose  $[f^{(k)}]^n$ and $[L_q(z,f)]^n$ share $(1,l)$ and $(\infty,m)$, where $0\leq l<\infty$ and $0\leq m\leq\infty$, and one of the following conditions holds:
	\begin{enumerate}
		\item [(i)] $l\geq 2$, $m=0$ and $n\geq 8$;
		\item [(ii)] $l\geq 2$, $m=\infty$ and $n\geq 7$;
		\item [(iii)] $l=1$, $m=0$ and $n\geq 9$;
		\item [(iv)] $l=0$, $m=0$ and $n\geq12$.
	\end{enumerate}
	Then $f^{(k)}=tL_q(z,f)$, for a non-zero constant $t$ that satisfies $t^n=1$.
\end{theo}
\par In 2018, Qi \textit{et al.} \cite{Qi & Li & Yang & 2018} also proved the following result.
\begin{theo}\label{t1.18}\cite{Qi & Li & Yang & 2018}
	Let $f$ be a meromorphic function of finite order. Suppose that $f^{\prime}$ and $\Delta_cf$ share $a_1, a_2, a_3, a_4$ IM, where $a_1, a_2, a_3, a_4$ are four distinct finite values. Then, $f^{\prime}(z)\equiv \Delta_cf.$
\end{theo}
\par We prove the following uniqueness theorem about the $k$-th derivative $f^{(k)}$ and linear difference polynomial $L(z,f)$ of a meromorphic function $f$, which is an extension of Theorem \ref{t1.18}.
\begin{theo}\label{t4}
	Let $f$ be a meromorphic function of finite order. Suppose that $f^{(k)}$ and $L(z,f)$ share $a_1, a_2, a_3, a_4$ IM, where $a_1, a_2, a_3, a_4$ are four distinct finite values. Then, $$f^{(k)}(z)\equiv L(z,f).$$
\end{theo}


\section{Key Lemmas}

In this section, we present some lemmas which will be needed in the sequel. Let $F$ and $G$ be two non-constant meromorphic functions defined in $\mathbb{C}$.  We also  denote by $H$, the following function \begin{eqnarray}\label{e2.1} H=\Big(\frac{F^{\prime\prime}}{F^{\prime}}-\frac{2F^{\prime}}{F-1}\Big)-\Big(\frac{G^{\prime\prime}}{G^{\prime}}-\frac{2G^{\prime}}{G-1}\Big).\end{eqnarray}
\begin{lem}\label{lem3.1}\cite{Lahiri & Complex Var & 2001}
	Let $F$, $G$ be two non-constant meromorphic functions such that they share $(1,1)$ and $H\not\equiv 0.$ Then \begin{eqnarray*} N(r,1;F\mid=1)=N(r,1;G\mid=1)\leq N(r,H)+S(r,F)+S(r,G).\end{eqnarray*}
\end{lem}
\begin{lem}\label{lem3.2}\cite{2}
	Let $F$, $G$ be two non-constant meromorphic functions sharing $(1,t),$ where $0\leq t<\infty.$ Then \begin{eqnarray*} && \ol N(r,1;F) + \ol N(r,1;G)-N_E^{1)}(r,1;F) +\left(t-\frac{1}{2}\right)\ol N_{*}(r,1;F,G) \\&\leq& \frac{1}{2} (N(r,1;F) + N(r,1;G)).\end{eqnarray*}
\end{lem}
\begin{lem}\label{lem3.3}\cite{Lahiri & Banerjee & 2006}
	Suppose $F$, $G$ share $(1,0)$, $(\infty,0)$. If $H\not \equiv 0,$ then, \beas  N(r,H)&&\leq N(r,0;F \mid\geq 2) + N(r,0;G \mid\geq 2)+\ol N_{*}(r,1;F,G) \\&&+ \ol N_{*}(r,\infty;F,G) + \ol N_{0}(r,0;F^{\prime}) + \ol N_{0}(r,0;G^{\prime})+S(r,F)+S(r,G),\eeas where $\ol N_{0}(r,0;F^{\prime})$ is the reduced counting function of those zeros of $F^{\prime}$ which are not the zeros of $F(F -1)$, and $\ol N_{0}(r,0;G^{\prime})$ is similarly defined. 
\end{lem}
\begin{lem}\label{lem3.4}\cite{Yang & 1972}
	Let $f$ be a non-constant meromorphic function and $P(f)=a_0+a_1f+a_2f^{2}+\ldots +a_{n}f^{n},$ where $a_0,a_1,a_2,\ldots, a_{n}$ are constants and $a_{n}\neq 0$. Then $$T(r,P(f))=nT(r,f)+O(1).$$
\end{lem}
\begin{lem}\label{lem3.5}\cite{Lahiri & Dewan & 2003} If $N\left(r,0;f^{(k)}\mid f\not=0\right)$ denotes the counting function of those zeros of  $f^{(k)}$ which are not the zeros of $f$, where a zero of $f^{(k)}$ is counted according to its multiplicity then $$N\left(r,0;f^{(k)}\mid f\not=0\right)\leq k\ol N(r,\infty;f)+N\left(r,0;f\mid <k\right)+k\ol N\left(r,0;f\mid\geq k\right)+S(r,f).$$\end{lem}
\begin{lem}\label{lem3.6}\cite{Yi & Kodai & 1999}
	Let $F$ and $G$ be two non-constant meromorphic functions such that they share $(1,0)$, and $H\not\equiv0$, then \beas N_E^{1)}(r,1;F)\leq N(r,\infty;H)+S(r,F)+S(r, G).\eeas Similar inequality holds for $G$ also. \end{lem}
\begin{lem}\label{lem3.7}\cite{Alzahare & Yi & 2004}
	If $F$, $G$ be two non-constant meromorphic functions such that they share $(1,1)$.
	Then \beas&& 2\ol N_L(r,1;F)+2\ol N_L(r,1;G)+\ol N_E^{(2}(r,1;F)-\ol N_{F>2}(r,1;G)\\&\leq& N(r,1;G)-\ol N(r,1;G).\eeas
\end{lem}
\begin{lem}\label{lem3.8}\cite{Banerjee & 2005}
	If two non-constant meromorphic functions $F$, $G$ share $(1,1)$, then 
	\beas \ol N_{F>2}(r,1;G)\leq \frac{1}{2}(\ol N(r,0;F)+\ol N(r,\infty;F)-N_{0}(r,0;F^{\prime}))+S(r,F),\eeas where $N_{0}(r,0;F^{\prime}))$ is the counting function of those zeros of $F^{\prime}$ which are not the zeros of $F(F-1)$.
\end{lem}
\begin{lem}\label{lem3.9}\cite{Banerjee & 2005}
	Let $F$ and $G$ be two non-constant meromorphic functions sharing $(1,0)$. Then \beas &&\ol N_L(r,1;F )+2\ol N_L(r,1;G)+\ol N_E^{(2}(r,1;F)-\ol N_{F>1}(r,1;G)-\ol N_{G>1}(r,1;F)\\&&\leq N(r,1;G)-\ol N(r,1;G).\eeas
\end{lem}
\begin{lem}\label{lem3.10}\cite{Banerjee & 2005}
	If $F$ and $G$ share $(1,0)$, then \beas \ol N_{L}(r,1;F)\leq \ol N(r,0;F)+\ol N(r,\infty;F)+S(r,F)\eeas
	\beas \ol N_{F>1}(r,1;G)\leq \ol N(r,0;F )+\ol N(r,\infty;F)-N_0(r,0;F^{\prime})+S(r,F).\eeas \text{Similar inequalities hold for} $G$ \text{also}.
\end{lem}
\begin{lem}\label{lem3.11}\cite{Yi & Kodai & 1999}
	Let $F$ and $G$ be two non-constant meromorphic functions such that they share $(1,0)$ and $H\not\equiv 0$. Then \beas N_E^{1)}(r,1;F)\leq N(r,\infty;H)+S(r,F)+S(r,H).\eeas
\end{lem}
\begin{lem}\label{lem3.12}\cite{Yang & Yi & Kluwer & 2003}
	Let $f$ and $g$ be two distinct non-constant rational
	functions and let $a_1, a_2, a_3, a_4$ be four distinct values. If $f$ and $g$ share $a_1, a_2, a_3, a_4$ IM, then $f(z)=g(z)$.
\end{lem}
\begin{lem}\label{lem3.13}\cite{Gundersen & Complex var & 1992}
	Suppose $f$ and $g$ are two distinct non-constant meromorphic functions, and $a_1, a_2, a_3, a_4\in\mathbb{C}\cup\{\infty\}$ are four distinct values. If $f$ and $g$ share $a_1, a_2, a_3, a_4$ IM, then
	\begin{enumerate}
		\item [(i)] $T(r,f)=T(r,g)+O(\log(rT(r,f)))$, as $r\not\in E$ and $r\rightarrow\infty$,
		\item [(ii)] $2T(r,f)=\sum_{j=1}^{4}\ol N\left(r,\displaystyle\frac{1}{f-a_j}\right)+O(\log(rT(r,f)))$, as $r\not\in E$ and $r\rightarrow\infty$, where
		$E\subset (1,\infty)$ is of finite linear measure.	
	\end{enumerate}
\end{lem}
\section{Proof of the theorems}
\par We prove only Theorems $\ref{t1}$ and $\ref{t4}$ as the proof of the rest of the theorems are very much similar to the proof of Theorem $\ref{t1}$.
\begin{proof}[\textbf{Proof of Theorem $\ref{t1}$}]
	\textbf{Case 1:} Suppose $H\not\equiv0$.\vspace{1mm} \par Let $F=(L(z,f))^n$ and $G=(f^{(k)})^n$.\vspace{1mm} 
	\par Keeping in view of Lemma $\ref{lem3.4}$, we get by applying Second fundamental theorem of Nevalinna on $F$ and $G$ that \bea\label{e4.1}&& n(T(r, L(z,f))+T(r,f^{(k)}))\nonumber\\&\leq& \ol N(r,0;F)+\ol N(r,1;F)+\ol N(r,\infty;F)+\ol N(r,0;G)+\ol N(r,1;G)\nonumber\\&&+\ol N(r,\infty;G)-\ol N_0(r,0;F^{\prime})-\ol N_0(r,0;G^{\prime})+S(r,F)+S(r,G),\eea where $\ol N_0(r,0;F^{\prime})$ and $\ol N_0(r,0;G^{\prime})$ are defined as in Lemma $\ref{lem3.3}$.\vspace{1mm}
	\par \textbf{(i).} Suppose $l\geq2$ and $m=0$.\vspace{1mm}
	\par Then using Lemmas $\ref{lem3.1}$, $\ref{lem3.2}$ and $\ref{lem3.3}$ in $(\ref{e4.1})$ we obtain
	\beas&& \frac{n}{2}(T(r,L(z,f))+T(r,(f^{(k)})))\\&\leq& N_2(r,0;F)+N_2(r,0;G)+\ol N(r,\infty;F)+\ol N(r,\infty;G)\\&&+\ol N_*(r,\infty;F,G)-\left(l-\frac{3}{2}\right)\ol N_{*}(r,1;F,G)+S(r,F)+S(r,G)\\&\leq& 2\ol N(r,0;L(z,f))+2\ol N(r,0; f^{(k)})+\ol N(r,\infty;L(z,f))\\&&+\ol N(r,\infty;f^{(k)})+\frac{1}{2}(\ol N(r,\infty;L(z,f))+\ol N(r,\infty;f^{(k)})\\&&+S(r,F)+S(r,G)\\&\leq&\frac{7}{2}(T(r,L(z,f))+T(r,f^{(k)})+S(r,F)+S(r,G).\eeas This implies that \beas (n-7)(T(r,L(z,f))+f^{(k)})\leq S(r,L(z,f))+S(r,f^{(k)}),\eeas which contradict to the fact that $n\geq8$.\vspace{1mm}
	\par \textbf{(ii).} Suppose $l\geq2$ and $m=\infty$.  Then using Lemmas $\ref{lem3.1}$, $\ref{lem3.2}$ and $\ref{lem3.3}$ in $(\ref{e4.1})$ we obtain
	\beas&& \frac{n}{2}(T(r,L(z,f))+T(r,f^{(k)}))\\&\leq& N_2(r,0;F)+N_2(r,0;G)+\ol N(r,\infty;F)+\ol N(r,\infty;G)\\&&-\left(l-\frac{3}{2}\right)\ol N_{*}(r,1;F,G)+S(r,F)+S(r,G)\\&\leq& 2\ol N(r,0;L(z,f))+2\ol N(r,0; f^{(k)})+\ol N(r,\infty;L(z,f))\\&&+\ol N(r,\infty;f^{(k)})+S(r,F)+S(r,G)\\&\leq& 3(T(r,L(z,f))+T(r,f^{(k)}))+S(r,F)+S(r,G).\eeas This implies that \beas (n-6)(T(r,L(z,f))+T(r,f^{(k)}))\leq S(r,L(z,f))+S(r,f^{(k)}),\eeas which contradict to the fact that $n\geq7$.\vspace{1mm}
	\par \textbf{(iii).} Suppose $l=1$ and $m=0$.
	\par Using Lemmas $\ref{lem3.1}$, $\ref{lem3.3}$, $\ref{lem3.7}$ and $\ref{lem3.8}$, we obtain
	\bea\label{e4.2}&& \ol N(r,1;F)\nonumber\\&\leq& N(r,1;F\mid=1)+\ol N_L(r,1;F)+\ol N_L(r,1;G)+\ol N_E^{(2}(r,1;F)\nonumber\\&\leq&\ol N(r,0;F\mid\geq2)+\ol N(r,0;G\mid\geq2)+\ol N_*(r,1;F,G)+\ol N_*(r,\infty;F,G)\nonumber\\&&+\ol N_L(r,1;F)+\ol N_L(r,1;G)+\ol N_E^{(2}(r,1;F)+\ol N_0(r,0;F^{\prime})+\ol N_0(r,0;G^{\prime})\nonumber\\&&+S(r,F)+S(r,G)\nonumber\\&\leq&\ol N(r,0;F\mid\geq2)+\ol N(r,0;G\mid\geq2)+\ol N_*(r,\infty;F,G)+2\ol N_L(r,1;F)\nonumber\\&&+2\ol N_L(r,1;G)+\ol N_E^{(2}(r,1;F)+\ol N_0(r,0;F^{\prime})+\ol N_0(r,0;G^{\prime})\nonumber\\&&+S(r,F)+S(r,G)\nonumber\\&\leq& \ol N(r,0;F\mid\geq2)+\ol N(r,0;G\mid\geq2)+\ol N_*(r,\infty;F,G)+N(r,1;G)\nonumber\\&&-\ol N(r,1;G)+\ol N_{F>2}(r,1;G)+\ol N_0(r,0;F^{\prime})+\ol N_0(r,0;G^{\prime})\nonumber\\&&+S(r,F)+S(r,G)\nonumber\\&\leq& \ol N(r,0;F\mid\geq2)+\ol N(r,0;G\mid\geq2)+\ol N_*(r,\infty;F,G)+N(r,0;G\mid G\neq0)\nonumber\\&&+\frac{1}{2}\ol N(r,0;F)+\frac{1}{2}\ol N(r,\infty;F)+\ol N_0(r,0;F^{\prime})+S(r,F)+S(r,G)\nonumber\\&\leq& \ol N(r,0;F\mid\geq2)+ N_2(r,0;G)+\ol N_*(r,\infty;F,G)+\ol N(r,\infty;G)\nonumber\\&&+\frac{1}{2}\ol N(r,0;F)+\frac{1}{2}\ol N(r,\infty;F)+\ol N_0(r,0;F^{\prime})+S(r,F)+S(r,G). \eea
	Similarly, we can get \bea\label{e4.3}&& \ol N(r,1;G)\nonumber\\&\leq&\ol N(r,0;G\mid\geq2)+ N_2(r,0;F)+\ol N_*(r,\infty;F,G)+\ol N(r,\infty;F)\nonumber\\&&+\frac{1}{2}\ol N(r,0;G)+\frac{1}{2}\ol N(r,\infty;G)+\ol N_0(r,0;G^{\prime})+S(r,F)+S(r,G).\eea
	Putting the values of $\ol N(r,1;F)$ and $\ol N(r,1;G)$ from $(\ref{e4.2})$ and $(\ref{e4.3})$ to $(\ref{e4.1})$, a simple calculation reduces to
	\beas&& n(T(r,L(z,f))+T(r,f^{(k)}))\\&\leq& 2N_2(r,0;F)+2N_2(r,0;G)+\frac{1}{2}(\ol N(r,0;F)+\ol N(r,0;G))\\&&+\frac{7}{2}(\ol N(r,\infty;F)+\ol N(r,\infty;G))+S(r,F)+S(r,G)\\&\leq& \frac{9}{2}(\ol N(r,0;L(z,f))+\ol N(r,0;f^{(k)})+\frac{7}{2}\ol N(r,\infty;L(z,f))\\&&+\frac{7}{2}\ol N(r,\infty;f^{(k)})+S(r,F)+S(r,G)\\&\leq& 8(T(r,L()z,f)+T(r,f^{(k)}))+S(r,L(z,f))+S(r,f^{(k)}),\eeas which is a contradiction since $n\geq9$.\vspace{1mm}
	\par \textbf{(iv).} Suppose $l=0$ and $m=0$.
	Using Lemmas $\ref{lem3.11}$, $\ref{lem3.3}$, $\ref{lem3.5}$, $\ref{lem3.9}$ and $\ref{lem3.10}$, we obtain 
	\bea\label{e4.4}&& \ol N(r,1;F)\nonumber\\&\leq& N_E^{1)}(r,1;F)+\ol N_L(r,1;F)+\ol N_L(r,1;G)+\ol N_E^{(2}(r,1;F)\nonumber\\&\leq& \ol N(r,0;F\mid\geq2)+\ol N(r,0;G\mid\geq2)+\ol N_*(r,1;F,G)+\ol N_*(r,\infty;F,G)\nonumber\\&&+\ol N_L(r,1;F)+\ol N_L(r,1;G)+\ol N_E^{(2}(r,1;F)+\ol N_0(r,0;F^{\prime})+\ol N_0(r,0;G^{\prime})\nonumber\\&&+S(r,F)+S(r,G)\nonumber\\&\leq& \ol N(r,0;F\mid\geq2)+\ol N(r,0;G\mid\geq2)+\ol N_*(r,\infty;F,G)+2\ol N_L(r,1;F)\nonumber\\&&+\ol 2N_L(r,1;G)+\ol N_E^{(2}(r,1;F)+\ol N_0(r,0;F^{\prime})+\ol N_0(r,0;G^{\prime})\nonumber\\&&+S(r,F)+S(r,G)\nonumber\\&\leq& \ol N(r,0;F\mid\geq2)+\ol N(r,0;G\mid\geq2)+\ol N_*(r,\infty;F,G)+\ol N_L(r,1;F)\nonumber\\&&+\ol N_{F>1}(r,1;G)+\ol N_{G>1}(r,1;F)+N(r,1;G)-\ol N(r,1;G)+\ol N_0(r,0;F^{\prime})\nonumber\\&&+\ol N_0(r,0;G^{\prime})+S(r,F)+S(r,G)\nonumber\\&\leq&
	\ol N(r,0;F\mid\geq2)+\ol N(r,0;G\mid\geq2)+\ol N_*(r,\infty;F,G)+2\ol N(r,0;F)\nonumber\\&&+2\ol N(r,\infty;F)+\ol N(r,0;G)+\ol N(r,\infty;G)+N(r,1;G)-\ol N(r,1;G)\nonumber\\&&+\ol N_0(r,0;F^{\prime})+\ol N_0(r,0;G^{\prime})+S(r,F)+S(r,G)
	\nonumber\\&\leq& N_2(r,0;F)+N_2(r,0;G)+\ol N(r,0;F)+2\ol N(r,\infty;F)+\ol N(r,\infty;G)\nonumber\\&&+\ol N_*(r,\infty;F,G)+N(r,0;G^{\prime}\mid G\neq0)+\ol N_0(r,0;F^{\prime})+S(r,F)+S(r,G)\nonumber\\&\leq&  N_2(r,0;F)+N_2(r,0;G)+\ol N(r,0;F)+\ol N(r,0;G)+2\ol N(r,\infty;F)\nonumber\\&&+2\ol N(r,\infty;G)+\ol N_*(r,\infty;F,G)+\ol N_0(r,0;F^{\prime})+S(r,F)+S(r,G). \eea
	Similarly, we can obtain 
	\bea\label{e4.5}&& \ol N(r.1;G)\nonumber\\&\leq& N_2(r,0;F)+N_2(r,0;G)+\ol N(r,0;F)+\ol N(r,0;G)+2\ol N(r,\infty;F)\nonumber\\&&+2\ol N(r,\infty;G)+\ol N_*(r,\infty;F,G)+\ol N_0(r,0;G^{\prime})+S(r,F)+S(r,G).\eea
	
	Using $(\ref{e4.4})$ and $(\ref{e4.4})$, $(\ref{e4.1})$ reduces to 
	\beas&& n(T(r,L(z,f))+T(r,f^{(k)}))\\&\leq& 2(N_2(r,0;F)+N_2(r,0;G))+3(\ol N(r,0;F)+\ol N(r,0;G))\\&&+2\ol N_*(r,\infty;F,G)+3(\ol N(r,\infty;F)+\ol N(r,\infty;G))\\&&+S(r,F)+S(r,G)\\&\leq& 7(\ol N(r,0;L(z,f))+\ol N(r,0;f^{(k)}))+4\ol N(r,\infty;L(z,f))\\&&+4\ol N(r,\infty;f^{(k)})+S(r,L(z,f))+S(r,f^{(k)})\\&\leq& 11(T(r,L(z,f))+T(r,f^{(k)}))+S(r,L(z,f))+S(r,f^{(k)}).\eeas
	This implies that 
	\beas(n-11)(T(r,L(z,f))+T(r,f^{(k)}))\leq S(r,L(z,f))+S(r,f^{(k)}),\eeas which is a contradiction since $n\geq12$.\vspace{1mm}
	\par \textbf{Case 2:} Suppose $H\equiv0$. Then by integration we get \bea \label{e4.6}F=\frac{AG+B}{CG+D},\eea where $A,\; B,\; C$ and $D$ are complex constants such that $AD-BC\neq0$.\vspace{1mm}
	\par From $(\ref{e4.6})$, it is easily seen that $T(r,L(z,f))=T(r,f^{(k)})+O(1)$.\vspace{1mm}
	\par \textbf{Subcase 2.1:} Suppose $AC\neq0$.
	Then $F-A/C=-(AD-BC)/C(CG+D)\neq 0.$ So $F$ omits the value $A/C.$ Therefore, by the second fundamental theorem, we get 
	\beas T(r,F)\leq \ol N(r,\infty;F)+\ol N(r,0;F)+\ol N\left(r,\frac{A}{C};F\right)+S(r,F).\eeas This implies that \beas nT(r,L(z,f))&\leq&\ol N(r,\infty; L(z,f))+\ol N(r,0;L(z,f))+S(r,L(z,f))\\&\leq& 2T(r,L(z,f))+S(r,L(z,f), \eeas which is not possible in all cases.\vspace{1mm}
	\par \textbf{Subcase 2.2:} Suppose that $AC=0$. Since $AD-BC\neq 0,$ both $A$ and $C$ can not be simultaneously zero.\vspace{1mm}
	\par \textbf{Subcase 2.2.1:} Suppose $ A\neq 0$ and $C=0.$ Then (\ref{e4.6}) becomes $ F \equiv \alpha G+\beta $, where $ \alpha=A/D$ and $\beta=B/D.$\vspace{1mm}
	\par If $F$ has no $1$-point, then by the second fundamental theorem of Nevalinna, we have \beas T(r,F)\leq \ol N(r,0;F)+\ol N(r,1;F)+\ol N(r,\infty;F)+S(r,F)\eeas $or,$ \beas (n-2)T(r,L(z,f))\leq S(r,L(z,f)),\eeas which is not possible in all cases.\vspace{1mm}
	\par Let $F$ has some $1$-point. Then $\alpha+\beta=1$. If $\beta=0,$ then $\alpha=1$ and then $F\equiv G$ which implies that $$L(z,f)=tf^{(k)},$$ where $t$ is a constant such that $t^n=1$.\vspace{1mm}
	\par Let $\beta\neq0$. Then applying the second main theorem of Nevalinna to $F$, we obtain
	\beas nT(r,L(z,f))&\leq& \ol N(r,\infty;F)+\ol N(r,0;F)+\ol N(r,\beta;F)+S(r,F)\\&\leq& 2T(r,L(z,f))+T(r,f^{(k)})+S(r,L(z,f))\\&\leq& 3T(r,L(z,f))+S(r,L(z,f)),\eeas which is not possible in all cases.\vspace{1mm}
	\par \textbf{Subcase 2.2.2:} Suppose $A=0$ and $C\neq 0$. Then (\ref{e4.6}) becomes \beas F\equiv \frac{1}{\gamma G_{1}+\delta},\eeas  where $\gamma=C/B$ and $\delta=D/B.$\par
	If $F$ has no $1$-point, then applying the second fundamental theorem to $F$, we have \beas nT(r,L(z,f))&\leq &\ol N(r,\infty;F)+\ol N(r,0;F)+\ol N(r,1;F)+S(r,F)\\&\leq & 2T(r,L(z,f))+ S(r,L(z,f)),\eeas which is a contradiction.\vspace{1mm}
	\par Suppose that $F$ has some $1$-point. Then $\gamma+\delta=1$.\par Therefore, $ F\equiv 1/(\gamma G+1-\gamma)$. Since $C\neq 0,$ $\gamma\neq 0$, and so $G$ omits the value $(\gamma-1)/\gamma.$\par
	By the second fundamental theorem of Nevalinna, we have \beas T(r,G)&\leq &\ol N(r,\infty;G)+\ol N(r,0;G)+\ol N\left(r,-\frac{1-\gamma}{\gamma};G\right)+S(r,G).\eeas $i.e.,$ \beas (n-2)T(r,f^{(k)})\leq S(r,f^{(k)}),\eeas which is a contradiction. This completes the proof of the theorem.\end{proof}
\begin{proof}[\textbf{Proof of Theorem $\ref{t4}$}]
	If $f$ is rational, the conclusion follows by Lemma \ref{lem3.12}. Assume that $f$ is transcendental meromorphic function. Then $f^{(k)}$ must be transcendental also. Now we discuss the following two cases.\vspace{1mm}
	\par \textbf{Case 1:} Suppose that $f^{(k)}$ is transcendental and $L(z,f)$ is rational. Then from Lemma \ref{lem3.13} (i), it follows that
	\beas T(r,f^{(k)})=T(r,L(z,f))+O(\log rT(r,f^{(k)}))=O(\log (rT(r,f^{(k)}))),\eeas which is a contradiction.\vspace{1mm}
	\par \textbf{Case 2:} Suppose $f$ and $L(z,f)$ are both transcendental.\vspace{1mm}
	\par Now keeping in view of Lemma \ref{lem3.13} (ii), and applying the second fundamental theorem of Nevalinna to $f^{(k)}$, we obtain
	\beas 3T(r,f^{(k)})&\leq& \ol N(r,f^{(k)})+\sum_{j=1}^{4}\ol N\left(r,\frac{1}{f^{(k)}-a_j}\right)+S(r,f^{(k)})\\&\leq& \ol N(r,f^{(k)})+2T(r,f^{(k)})+S(r,f^{(k)})\\&\leq& 3T(r,f^{(k)})+S(r,f^{(k)}),\eeas which implies that 
	\beas N(r,f^{(k)})=\ol N(r,f^{(k)})+S(r,f^{(k)}).\eeas
	This implies that \beas N(r,f)+k\ol N(r,f)=N(r,f^{(k)})=\ol N(r,f^{(k)})+S(r,f^{(k)})=\ol N(r,f)+S(r,f^{(k)}).\eeas This shows that \bea\label{e4.7} N(r,f)=\ol N(r,f)=\ol N(r,f^{(k)})=S(r,f^{(k)}).\eea Again from Lemma \ref{lem3.13} (i), we have 
	\beas T(r,f^{(k)})=T(r,L(z,f))+S(r,f^{(k)}).\eeas
	Keeping in view of (\ref{e4.7}), the above equation,  and applying the second main theorem to $f^{(k)}$, we obtain
	\beas 3T(r,f^{(k)})&\leq& \ol N(r,f^{(k)})+\sum_{j=1}^{4}\ol N\left(r,\frac{1}{f^{(k)}-1}\right)+S(r,f^{(k)})\\&\leq& \ol N(r,f^{(k)})+N\left(r,\frac{1}{f^{(k)}-L(z,f)}\right)+S(r,f^{(k)})\\&\leq& T(r,f^{(k)})+T(r,L(z,f))+S(r,f^{(k)})\\&\leq& 2T(r,f^{(k)})+S(r,f^{(k)}),\eeas which is a contradiction. This completes the proof of the theorem.	
\end{proof}
\bigskip


\begin{thebibliography}{99}
\bibitem{Banerjee & Bhattacharya & 2019} A. Banerjee and S. Bhattacharayya, On the uniqueness of meromorphic functions and its difference operators sharing values of sets, Rend. Circ. Mat. Palermo2. Ser., doi 10.1007/s12215-016-0295-1.
\end{thebibliography}
\end{document}